\documentclass{agtart_a}
\pdfoutput=1
\usepackage{pinlabel}

\title{Sweepouts of amalgamated 3--manifolds}

\author{David Bachman}
\givenname{David}
\surname{Bachman}
\address{Mathematics Department\\
Pitzer College\\\newline
1050 North Mills Avenue\\
Claremont CA 91711\\USA}
\email{bachman@pitzer.edu}

\author{Saul Schleimer}
\givenname{Saul}
\surname{Schleimer}
\address{Department of Mathematics\\Rutgers\\
The State University of New Jersey\\\newline
110 Frelinghuysen Rd\\
Piscataway NJ 08854-8019\\USA}
\email{saulsch@math.rutgers.edu}

\author{Eric Sedgwick}
\givenname{Eric}
\surname{Sedgwick}
\address{CTI\\DePaul University\\
243 S Wabash Avenue\\\newline
Chicago IL 60604\\USA}
\email{esedgwick@cs.depaul.edu}



\received{26 July  2005}
\revised{18 January 2006}
\accepted{26 January 2006}
\proposed{}
\seconded{}
\publishedonline{24 February 2006}
\published{24 February 2006}

\volumenumber{6}
\issuenumber{}
\publicationyear{2006}
\papernumber{6}
\startpage{171}
\endpage{194}

\doi{}
\MR{}
\Zbl{}

\newcommand{\calF}{{\mathcal{F}}}

\newcommand{\TT}{{\mathbb{T}}}
\newcommand{\VV}{{\mathbb{V}}}
\newcommand{\WW}{{\mathbb{W}}}
\newcommand{\ZZ}{{\mathbb{Z}}}

\renewcommand{\setminus}{{\smallsetminus}}

\let\from\co   

\newcommand{\homeo}{{\medspace \cong \medspace}} 
 
\newcommand{\isom}{{\medspace \cong \medspace}} 

\newcommand{\bdy}{{\partial}} 

\newcommand{\MCG}{{\mathcal{MCG}}} 

\makeatletter
\def\cnewtheorem#1[#2]#3{\newtheorem{#1}{#3}[section]
\expandafter\let\csname c@#1\endcsname\c@theorem}

\theoremstyle{plain}
\newtheorem{theorem}{Theorem}[section]
\cnewtheorem{corollary}[theorem]{Corollary}
\cnewtheorem{lemma}[theorem]{Lemma}
\cnewtheorem{conjecture}[theorem]{Conjecture}

\theoremstyle{definition}
\cnewtheorem{define}[theorem]{Definition}
\cnewtheorem{claim}[theorem]{Claim}

\newtheorem*{proofclaim}{Claim}
\cnewtheorem{proposition}[theorem]{Proposition}
\cnewtheorem{remark}[theorem]{Remark}
\cnewtheorem{remarks}[theorem]{Remarks}

\makeatother  

\newsavebox{\savepar}

\newcommand{\band}[1]{{\widehat{#1}}} 

\begin{document}

\begin{asciiabstract}
We show that if two 3-manifolds with toroidal boundary are glued via a
`sufficiently complicated' map then every Heegaard splitting of the
resulting 3-manifold is weakly reducible.  Additionally, if Z is a
manifold obtained by gluing X and Y, two connected small manifolds
with incompressible boundary, along a closed surface F. Then the genus
g(Z) of Z is greater than or equal to 1/2(g(X)+g(Y)-2g(F)).  Both
results follow from a new technique to simplify the intersection
between an incompressible surface and a strongly irreducible Heegaard
splitting.
\end{asciiabstract}

\begin{mathmlabstract}
We show that if two 3&ndash;manifolds with toroidal boundary are glued via a
"sufficiently complicated" map then every Heegaard splitting of the
resulting 3&ndash;manifold is weakly reducible. Additionally, suppose
<math overflow="scroll">
    <mi>X</mi>
    <msub>
      <mo>&#x222A;</mo>
      <mi>F</mi>
    </msub>
    <mi>Y</mi>
</math>
is a manifold obtained by gluing
<math overflow="scroll">
    <mi>X</mi>
</math>
and
<math overflow="scroll">
    <mi>Y</mi>
</math>, two connected
small manifolds with incompressible boundary, along a closed surface
<math overflow="scroll">
    <mi>F</mi>
</math>. Then the following inequality on genera is obtained:

  <math display="block" overflow="scroll">
    <mi>g</mi>
    <mo>&#x2061;</mo>
    <mfenced separators="">
      <mi>X</mi>
      <msub>
        <mo>&#x222A;</mo>
        <mi>F</mi>
      </msub>
      <mi>Y</mi>
    </mfenced>
    <mo>&#x2265;</mo>
    <mfrac>
      <mn>1</mn>
      <mn>2</mn>
    </mfrac>
    <mfenced close=")" open="(" separators="">
      <mi>g</mi>
      <mo>&#x2061;</mo>
      <mfenced>
        <mi>X</mi>
      </mfenced>
      <mo>+</mo>
      <mi>g</mi>
      <mo>&#x2061;</mo>
      <mfenced>
        <mi>Y</mi>
      </mfenced>
      <mo>-</mo>
      <mrow>
        <mn>2</mn>
        <mi>g</mi>
      </mrow>
      <mfenced>
        <mi>F</mi>
      </mfenced>
    </mfenced>
    <mi>.</mi>
</math>
Both results follow from a new technique to simplify the 
intersection between an incompressible surface and a 
strongly irreducible Heegaard splitting.
\end{mathmlabstract}

\begin{abstract}
We show that if two 3--manifolds with toroidal boundary are glued via a
``sufficiently complicated" map then every Heegaard splitting of the
resulting 3--manifold is weakly reducible. Additionally, suppose $X
\cup _F Y$ is a manifold obtained by gluing $X$ and $Y$, two connected
small manifolds with incompressible boundary, along a closed surface
$F$. Then the following inequality on genera is obtained:
\[g(X \cup _F Y)\ge \tfrac{1}{2}\left(g(X)+g(Y)-2g(F)\right).\]
Both results follow from a new technique to simplify the intersection between an incompressible surface and a strongly irreducible Heegaard splitting.
\end{abstract}

\maketitle

\section{Introduction}

It is a consequence of the Haken Lemma \cite{haken:68} and the Uniqueness of Prime Decompositions, Kneser \cite{kneser:29}, that Heegaard genus is well behaved under connected sum.  In particular, 3--manifold genus is additive: 
\[g(X \# Y) = g(X) + g(Y)\]
Here we discuss the Heegaard splittings of a manifold obtained by gluing together manifolds along boundary components of higher genus.

To this end let $X$ and $Y$ be 3--manifolds with incompressible boundary homeomorphic to a connected surface $F$. It is not difficult to show that if $H_X$ and $H_Y$ are Heegaard surfaces in $X$ and $Y$ then we can amalgamate these splittings to obtain a Heegaard surface in $X \cup _F Y$ with genus equal to $g(H_X)+g(H_Y)-g(F)$ (see, for example, Schultens \cite{schultens:93}). Letting $g(X)$, $g(Y)$, and $g(X \cup_F Y)$ denote the minimal genus among all Heegaard surfaces in the respective 3--manifolds, we find:
\begin{equation}
\label{e:ineq}
g(X \cup_F Y)\le g(X)+g(Y)-g(F)
\end{equation}
Bounds in the other direction are harder to obtain. When $F \cong S^2$ it follows from the Haken Lemma \cite{haken:68} that the above inequality may be replaced by an equality.  In \fullref{s:saul} we examine the case where $F$ is a torus.   We assume here that the map which identifies  $\partial X$ to $\partial Y$ is ``sufficiently complicated," in a sense to be made precise in \fullref{s:saul}. 

\medskip
\noindent {\bf  
\fullref{Thm:ManifoldsWithoutStronglyIrreducibleSplittings}}\qua {\sl
Suppose that $X$ and $Y$ are knot manifolds and $\varphi \from \bdy X
\to \bdy Y$ is a sufficiently complicated homeomorphism.  Then the manifold
$M(\varphi) = X \cup_\varphi Y$ has no strongly irreducible
Heegaard splittings.}
\medskip

\noindent In particular it follows from this result that every Heegaard splitting of $X \cup_F Y$ is an amalgamation of splittings of $X$ and $Y$. In this situation Inequality \ref{e:ineq} becomes an equality.

In the case where the genus of $F$ is at least two there is the following result of Lackenby ~\cite{lackenby:03}:

\medskip
\noindent {\bf Theorem}\qua {\sl
Let $X$ and $Y$ be simple 3--manifolds, and let $h\co \partial X \to F$
and $h'\co F \to \partial Y$ be homeomorphisms with some connected
surface $F$ of genus at least two. Let $\psi\co  F \to F$ be a
psuedo-Anosov homeomorphism. Then, provided $|n|$ is sufficiently
large,
\[g(X \cup_{h'\psi^n h} Y)=g(X)+g(Y)-g(F).\]
Furthermore, any minimal genus Heegaard splitting for $X
\cup_{h'\psi^n h} Y$ is obtained from splittings of $X$ and $Y$ by
amalgamation, and hence is weakly reducible.}
\medskip

If $\psi$ fails to be ``sufficiently complicated" then
there is no hope of an exact equality, as in the previous
theorem. 
Previous known lower bounds were obtained by
Johannson \cite{johannson:95} when $X$ and $Y$ are simple
\[g(X \cup_F Y)\ge \frac{1}{5}(g(X)+g(Y)-2g(F)).\]
Schultens has generalized this result to allow essential annuli \cite{schultens:01}.

By assuming the component manifolds $X$ and $Y$ are {\it small} we get a new bound. The following statement is one case of \fullref{Cor:genusbounds}:

\medskip
\noindent {\bf \fullref{Cor:genusbounds}$'$}\qua {\sl
Suppose $X$ and $Y$ are compact, orientable, connected, small 3--manif\-olds with incompressible boundary homeomorphic to a surface $F$. Then
\[g(X \cup_F Y)\ge \frac{1}{2}(g(X)+g(Y)-2g(F)).\]}
\eject

Both of our results follow from showing that a strongly
irreducible Heegaard surface $H$ can be isotoped to meet
the gluing surface $F$ in a particularly nice fashion. Often in
these types of arguments one simplifies the intersection by making
every loop of $H \cap F$ essential in both surfaces. In this
paper, rather than focusing on the intersection set $H \cap F$, we
focus on the complimentary pieces $H \setminus N(F)$. Our result is that
$H$ and $F$ may always be arranged so that {\it almost} every component $H'$ of $H
\setminus N(F)$ is {\it incompressible}. On such a component every loop which is
essential in $H'$ is essential in $M \setminus N(F)$. There is at most one component $H''$ which is compressible. In this case we find that $H''$ is {\it
strongly irreducible}, in the sense that every
essential loop which bounds a disk on one side meets every
essential loop bounding a disk on the other. See
\fullref{Lem:StronglyIrreducibleSweepout}.

\section{Definitions}
In this section we give some of the standard definitions that will be used throughout paper. 

\subsection{Essential loops, arcs, and surfaces}
A loop $\gamma$ embedded in the interior of a compact, orientable surface $F$ is called
{\it essential} if it does not bound a disk in $F$. If $F$ is embedded
in a 3--manifold, $M$, a {\it compressing disk} for $F$ is a disk, $D
\subset M$, such that $F \cap D=\partial D$, and such that $\partial
D$ is essential on $F$. If we identify a thickening of $D$ in
$M \setminus N(F)$ with $D \times I$ then to {\it compress $F$
along $D$} is to remove $(\partial D) \times I$ from $F$ and replace
it with $D \times \partial I$.

A properly embedded arc $\alpha$ on $F$ is {\it essential} if there is no subarc $\beta$ of $\partial F$ such that $\alpha \cup \beta$ is the boundary of a subdisk of $F$. If $F$ is properly embedded in a 3--manifold, $M$, a {\it boundary-compressing disk} is a disk, $D$, such that $\partial D=\alpha \cup \beta$, where $F \cap D=\alpha$ is an essential arc on $F$ and $D \cap \partial M=\beta$. If we identify a thickening of $D$ in $M \setminus N(F)$ with $D \times I$ then to {\it boundary-compress $F$
along $D$} is to remove $\alpha \times I$ from $F$ and replace
it with $D \times \partial I$.

A properly embedded surface is {\it incompressible} if there are no
compressing disks for it. A properly embedded, separating surface is
{\it strongly irreducible} if there are compressing disks for it on
both sides, and each compressing disk on one side meets each
compressing disk on the other side.

A compact, orientable 3--manifold is said to be {\it irreducible} if every embedded
2--sphere bounds a 3--ball. A 3--manifold is said to be {\it small} if it
is irreducible and every incompressible surface is parallel to a
boundary component.

\subsection{Heegaard and generalized Heegaard Splittings.}
A {\it compression body} is a 3--manifold $C$ constructed in one of two different ways. The first way is to begin with a collection of zero--handles and attach one--handles to their boundaries, resulting in a manifold that may or may not be connected. In this case we say the {\it spine} of $C$ is a 1--complex $\Sigma$ in $C$ such that $C$ is homeomorphic to a thickening of $\Sigma$. We set $\partial _-C=\emptyset$ and $\partial _+C=\partial C$. 

The second way to construct a compression body is to begin with a closed (possibly disconnected) orientable surface $F$ with no sphere components, and let $C$ be the manifold obtained by attaching one--handles to the surface $F \times \{1\} \subset F \times I$. In this case we say $\partial_-C = F \times \{0\}$ and $\partial_+ C =\partial C  \setminus \partial_-C$. The spine $\Sigma$ is then the union of $\partial _-C$ and a collection of arcs which are properly embedded in $C$, such that $C$ is a thickening of $\Sigma$. 

A surface, $H$, in a 3--manifold, $M$, is a {\it Heegaard surface for
  M} if $H$ separates $M$ into two compression bodies, $V$ and $W$,
such that $H=\partial _+V=\partial _+ W$.

A {\it generalized Heegaard splitting} of a 3--manifold $M$, Scharlemann--Thompson
\cite{ScharlemannThompson94}, is a sequence
$\{H_i\}_{i=0}^{2n}$ of pairwise disjoint, closed surfaces in $M$
such that
\begin{itemize}
\item $\partial M=H_0 \cup H_{2n}$ (if $\partial M=\emptyset$ then
  $H_0=H_{2n}=\emptyset$) and 
\item for each odd $i$, the surface $H_i$ is a Heegaard splitting of the submanifold co\-bounded by $H_{i-1}$ and $H_{i+1}$.
\end{itemize}
We will call the
set of surfaces with even index {\it thin levels} and the set with odd
index {\it thick levels}. 

Generalized Heegaard splittings are associated to handle structures in the following way. Given a generalized Heegaard splitting $\{H_i\}_{i=0}^{n}$ there is a sequence of submanifolds $\{M_i\}$ of $M$ as follows:

\begin{itemize}
	\item $M_0$ is a union of zero--handles and 1--handles.
	\item For odd $i$ between $1$ and $n$, $M_i$ is obtained from $M_{i-1}$ by attaching one--handles.
	\item For even $i$ between $2$ and $n-1$, $M_i$ is obtained from $M_{i-1}$ by attaching two--handles.
	\item $M_n=M$ is obtained from $M_{n-1}$ by attaching two--handles and three--handles.
\end{itemize}
Conversely, given a handle structure for $M$ there is an associated generalized Heegaard splitting as above. 

Suppose $H_X$ and $H_Y$ are Heegaard surfaces in 3--manifolds $X$ and $Y$. Suppose further that the boundaries of both $X$ and $Y$ are homeomorphic to a surface $F$. Then $\{\emptyset, H_X, F, H_Y, \emptyset\}$ is a generalized Heegaard splitting of $X \cup _F Y$. We may now choose a handle structure associated to this generalized Heegaard splitting, and re-arrange it so that handles are added in order of increasing index. The generalized Heegaard splitting associated to this new handle structure will be of the form $\{\emptyset, H, \emptyset\}$, where $H$ is a Heegaard surface in $X \cup _F Y$. In this case the Heegaard surface $H$ is the {\it amalgamation} of $H_X$ and $H_Y$, as defined by Schultens \cite{schultens:93}. 

\subsection{Normal and almost normal surfaces.}
\begin{figure}[ht!]
        \begin{center}
        \includegraphics[width=3in]{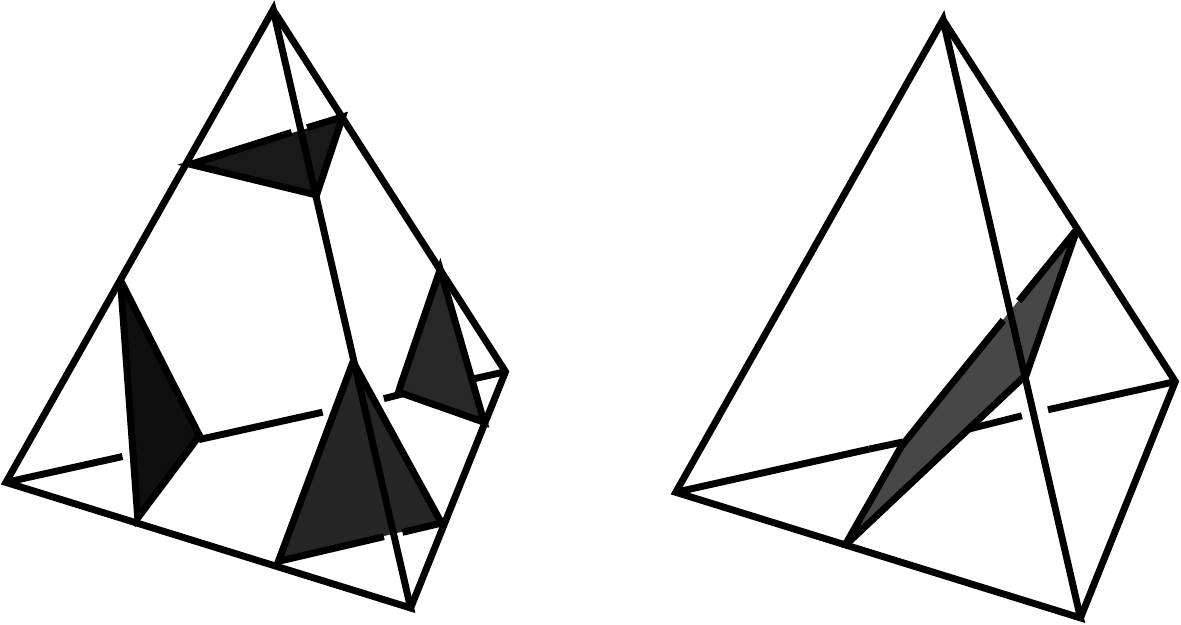}
        \caption{Normal disks}
        \label{Fig:normal}
        \end{center}
        \end{figure}

A {\it normal disk} in a tetrahedron is a triangle or a quadrilateral, as in \fullref{Fig:normal}.  Let $X$ be a 3--manifold equipped with a psuedo-triangulation. That is, $X$ is expressed as a collection of tetrahedra, together with face pairings. 

A properly embedded surface in $X$ is {\it normal} if it intersects every tetrahedron in a collection of {\it triangles} and {\it quadrilaterals}. Normal surfaces were first introduced by Kneser \cite{kneser:29}, and later used to solve several important problems by Haken \cite{haken:61}.

      \begin{figure}[ht!]
        \vspace{0 in}
        \cl{%
        \psfig{file=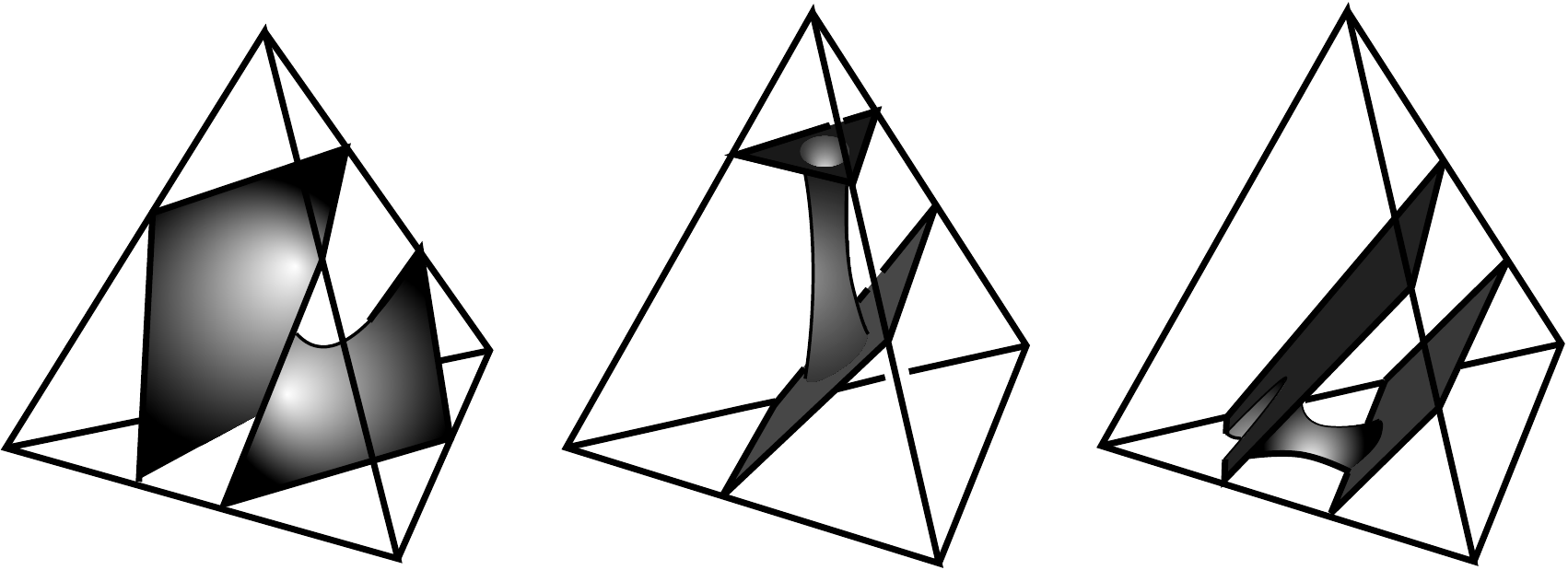,width=4in}}
        \caption{Exceptional disks in an almost normal surface}
        \label{Fig:almostnormal}
        \end{figure}

A properly embedded surface in $X$ is {\it almost normal} if it is normal everywhere, with the exception of exactly one piece in one tetrahedron. The exceptional piece can either be an octagon, two normal disks connected by an unknotted tube, or two normal disks connected by a band along $\partial X$ (see \fullref{Fig:almostnormal}). In the closed case, almost normal surfaces were introduced by Rubinstein \cite{rubinstein:95}. They were later generalized to surfaces with non-empty boundary by the first author \cite{bachman:01}.

\section{Labelling sweepouts}
\label{s:labelling}

In this section we prove the technical lemmas on which Sections \ref{s:saul} and \ref{s:eric} rely.

\begin{lemma}[Scharlemann \cite{scharlemann:98}]
\label{Lem:NoNesting}
Let $H$ be a strongly irreducible Heegaard surface, and $\gamma$ be an
essential curve on $H$. Suppose $\gamma$ bounds a disk $D \subset M$ such
that $D$ is transverse to $H$. Then $\gamma$ bounds a compressing disk
for $H$.
\end{lemma}


\begin{define}
Two surfaces $H$ and $F$ embedded in a 3--manifold are {\it almost
transverse} if they have exactly one non-transverse intersection point, and it is a
saddle point.
\end{define}

\begin{lemma} 
\label{Lem:StronglyIrreducibleSweepout}
Let $M$ be a compact, irreducible, orientable 3--manifold with $\partial M$ incompressible, if non-empty. Suppose $M=V \cup _H W$,
where $H$ is a strongly irreducible Heegaard surface. Suppose further
that $M$ contains an incompressible, orientable, closed, non-boundary
parallel surface $F$. Then either
\begin{itemize}
	\item $H$ may be isotoped to be transverse
to $F$, with every component of $H \setminus N(F)$ incompressible in the
respective submanifold of $M \setminus N(F)$,
	\item $H$ may be isotoped to be transverse
to $F$, with every component of $H \setminus N(F)$ incompressible in the
respective submanifold of $M \setminus N(F)$ except for exactly one
strongly irreducible component, or
	\item $H$ may be isotoped to be almost transverse
to $F$, with every component of $H \setminus N(F)$ incompressible in the
respective submanifold of $M \setminus N(F)$.
\end{itemize}
\end{lemma}

\begin{remarks}\quad
\label{Rmk:HTransverseToF}
\begin{enumerate}
\item After applying the lemma every loop of $H \cap F$ must be essential on both surfaces. Otherwise there is such a loop that is inessential on $F$ and essential on $H$. This loop, after a small isotopy, bounds a compressing disk $D$ for a component $H'$ of $H \setminus N(F)$. By the lemma, $H'$ must then be strongly irreducible. But $D$ is disjoint from every compressing disk for $H'$ on the opposite side, a contradiction.
\item In the case where $F \cong \TT^2$ it will follow from the proof
  that $H$ may actually be isotoped to be transverse to $F$. Here, only conclusions one or two of the lemma occur. 
\end{enumerate}
\end{remarks}

\begin{proof}[Proof of \fullref{Lem:StronglyIrreducibleSweepout}]
Choose spines $\Sigma _V$ of $V$ and $\Sigma _W$ of $W$.

\begin{claim}
\label{Clm:FMeetsSpines}
The surface $F$ meets both $\Sigma _V$ and $\Sigma _W$.
\end{claim}

\begin{proof}
Suppose $F \cap \Sigma _V=\emptyset$. Then $F$ lies in a compression
body homeomorphic to $W$. As the only incompressible surfaces in $W$
are components of $\partial _-W$, we conclude that $F$ is boundary
parallel in $M$. This violates the hypotheses of
\fullref{Lem:StronglyIrreducibleSweepout}.
\end{proof}

Fix a {\it sweepout} of $M$: a continuous map $\Phi \co  H \times
I \rightarrow M$ such that
\begin{itemize}
\item $H(0)=\Sigma _V$,
\item $H(1)=\Sigma _W$, and
\item the restriction of $\Phi$ to $H \times (0,1)$ is a
  smooth homeomorphism onto the complement of $\Sigma _V \cup \Sigma
  _W$.
\end{itemize}
Here $H(t)=\Phi(H \times t)$.
The map $\Phi$ is a {\it sweepout} of $M$. (Note that this is a slightly different definition than the one introduced by Rubinstein). Let $V(t)$ and $W(t)$
denote the compression bodies bounded by $H(t)$ (where $\Sigma _V
\subset V(t)$).


The sweepout $\Phi$ induces a height function $h\co  F \to I$ as follows. Define $h(x)=t$ if $x \in \Phi(H,t)$. Perturb $F$ so that $h$ is Morse on $F \setminus (\Sigma _V \cup \Sigma _W)$. Let $\{t_i\}_{i=0}^n$ denote the set of critical values of $h$. It follows from
\fullref{Clm:FMeetsSpines} that $t_0=0$ and $t_n=1$. 
We now label
each subinterval $(t_i,t_{i+1})$ with the letters $\VV$ and/or $\WW$
by the following scheme. If, for some $t \in (t_i,t_{i+1})$, there is
a compressing disk for $H(t)$ in $V(t)$ with boundary disjoint
from $F$ then label this subinterval with the letter $\VV$. See \fullref{Fig:labelling}. Similarly,
if there is a compressing disk in $W(t)$ with boundary disjoint
from $F$ then label with the letter $\WW$.

        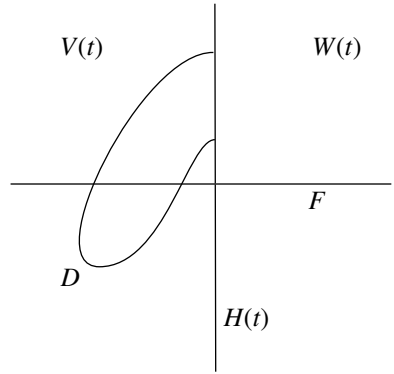
\begin{figure}[ht!]
        \labellist\small
        \pinlabel $H(t)$ [l] at 308 411
        \pinlabel $V(t)$ at 206 617
        \pinlabel $W(t)$ at 400 617
        \pinlabel $F$ [t] at 385 514
        \pinlabel $D$ [tr] at 209 454
        \endlabellist        
        \cl{%
        \includegraphics[width=2in]{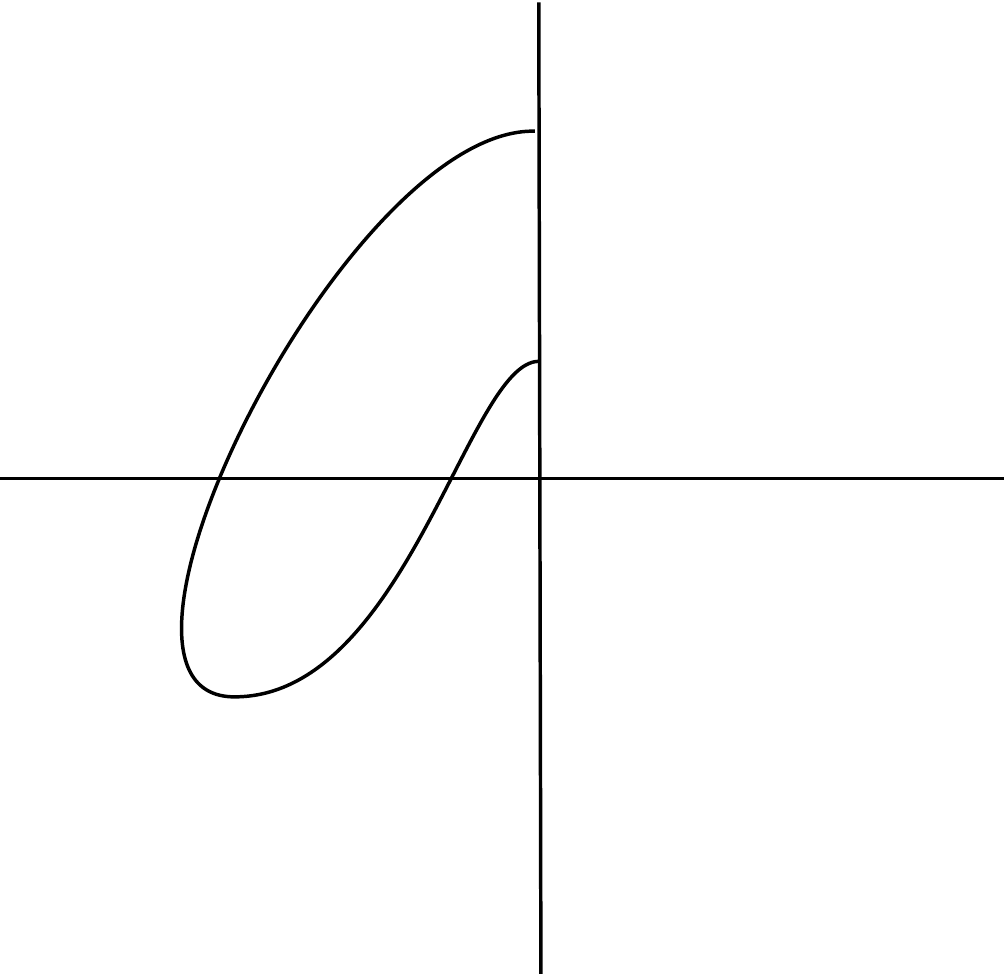}}
        \caption{If $D$ is a compressing disk for $H(t)$ in $V(t)$ with boundary disjoint from $F$ then the interval containing $t$ would get the label $\VV$.}
        \label{Fig:labelling}
        \end{figure}

\begin{claim}
\label{Clm:UnlabelledInterval}
If the subinterval $(t_i,t_{i+1})$ is unlabelled then the first conclusion
of \fullref{Lem:StronglyIrreducibleSweepout} follows.
\end{claim}

\begin{proof}
Suppose $t \in (t_i,t_{i+1})$. First, we claim that all curves of
$H(t) \cap F$ are essential on both or inessential on both. If not
then, as $F$ is incompressible, there is a loop $\delta \subset H(t)
\cap F$ that is inessential on $F$ but essential on $H(t)$. The loop
$\delta$ bounds a disk $D \subset F$.  Thus the hypotheses of
\fullref{Lem:NoNesting} are satisfied.  It follows that $\delta$
bounds a compressing disk in $V(t)$ or in $W(t)$.  Finally, $\delta$
may be isotoped inside of $H(t)$  by a small pushout move to be disjoint from $F$.  This
violates the assumption that $(t_i,t_{i+1})$ is unlabelled. We deduce that all curves of $H(t) \cap F$ are essential or inessential on both. 


As $M$ is irreducible we may isotope $H(t)$ to remove those loops of
$H(t) \cap F$ which are inessential on both surfaces, without
affecting those loops of $H(t) \cap F$ which were essential on both.
We now claim that after such an isotopy any essential loop of $H(t)
\setminus N(F)$ is essential on $H(t)$. We prove the contrapositive:
Suppose $E \subset H(t)$ is an embedded disk with $\bdy E \cap F =
\emptyset$.  All curves of $E \cap F$ are inessential on both
surfaces. Isotope $E$ rel boundary to make $E \cap F = \emptyset$. We conclude $E \subset M
\setminus N(F)$, and hence $\bdy E$ is inessential on $H(t) \setminus N(F)$.

Finally, we claim that the components of $H(t) \setminus N(F)$ are
incompressible in the respective submanifolds of $M \setminus N(F)$.
Suppose $H'$ is a compressible component. Then there is an essential
loop $\gamma \subset H'$ which bounds a compressing disk for $H'$.  By
the preceding remarks $\gamma$ is essential on $H(t)$ as well. By
\fullref{Lem:NoNesting} the loop $\gamma$ bounds a compressing disk
for $H(t)$, which must be in $V(t)$ or $W(t)$. This now contradicts
the fact that $(t_i,t_{i+1})$ is unlabelled.
\end{proof}

\begin{claim}
\label{Clm:BothLabels}
If the subinterval $(t_i,t_{i+1})$ has both of the labels $\VV$ and
$\WW$ then the second conclusion of
\fullref{Lem:StronglyIrreducibleSweepout} follows.
\end{claim}

\begin{proof}
Suppose $t \in (t_i,t_{i+1})$. We begin as in the proof of
\fullref{Clm:UnlabelledInterval} by asserting that all curves of
$H(t) \cap F$ are either inessential or essential on both. If not,
then as above there is a loop $\delta \subset H(t) \cap F$ which
bounds a compressing disk for $H(t)$. Suppose $\delta$ bounds a
compressing disk in $V(t)$. (The other case is similar.) Since $(t_i,t_{i+1})$ has the label $\WW$
there is a loop $\gamma$ on some component of $H(t) \setminus N(F)$ which
bounds a disk in $W(t)$. But then $\delta \cap \gamma =\emptyset$
contradicts the strong irreducibility of $H$.

As in the proof of \fullref{Clm:UnlabelledInterval} it now follows
that we may isotope $H(t)$, preserving the set of loops of $H(t) \cap
F$ which are essential on both, so that any loop which is essential on
$H(t) \setminus N(F)$ is also essential on $H(t)$.

Let $H'$ be a component of $H(t) \setminus N(F)$ which contains a loop
$\gamma$ bounding a compressing disk for $H(t)$ in $W(t)$.  By strong
irreducibility of $H(t)$ any essential loop of $H(t) \setminus N(F)$ which bounds a
compressing disk in $V(t)$ must meet $\gamma$, and hence must also lie
in $H'$. Furthermore, since the subinterval $(t_i,t_{i+1})$ has the
label $\VV$, there is at least one such loop $\rho$. By identical
reasoning we conclude that any essential loop of $H(t) \setminus N(F)$
which bounds a compressing disk in $W(t)$ must meet $\rho$, and hence
must also be on $H'$. We conclude that there are no loops on any other
component of $H(t) \setminus N(F)$ which bound compressing disks.  Hence
all components of $(H(t) \setminus N(F)) \setminus H'$ are incompressible
in the respective submanifolds of $M \setminus N(F)$.  Furthermore, the
strong irreducibility of $H'$ follows from the existence of the $\VV$
and $\WW$ labels and strong irreducibility of $H(t)$.
\end{proof}

\begin{claim}
\label{Clm:LabelChangeImpliesSaddle}
If the labelling of $(t_{i-1}, t_i)$ is different from that of
$(t_i,t_{i+1})$ then the critical value $t_i$ corresponds to a saddle
tangency between $H(t_i)$ and $F$. \qed
\end{claim}

\begin{claim}
\label{Clm:BeginVEndW}
The subinterval $(0,t_1)$ is labelled $\VV$ and the subinterval
$(t_{n-1},1)$ is labelled $\WW$.
\end{claim}

\begin{proof}
For sufficiently small $\epsilon$ the surface $H(\epsilon)$ looks
like the frontier of a neighborhood of $\Sigma _V$. By
\fullref{Clm:FMeetsSpines} the surface $F$ meets $\Sigma _V$.  Hence, $F$ contains small compressions for
$H(\epsilon)$ in $V(\epsilon)$. We can push these compressions off
$F$, giving compressions with boundary on a component of
$H(\epsilon) \setminus N(F)$ in $V(\epsilon)$. Hence, the label of
$(0,t_1)$ is $\VV$. A symmetric argument completes the proof of the
claim.
\end{proof}

Following Claims~\ref{Clm:UnlabelledInterval} and~\ref{Clm:BothLabels}
we now assume that every subinterval has a label. Furthermore, we assume that every subinterval has exactly one label: either $\VV$ or $\WW$, but not both.
It then follows from \fullref{Clm:BeginVEndW} that there is
some first critical value $t_i$ where the labelling changes from $\VV$ to
$\WW$. By \fullref{Clm:LabelChangeImpliesSaddle} this critical value
must correspond to a saddle tangency. 

\begin{claim}
There is a surface $H_0$, isotopic to $H(t_i)$, such that all components of $H_0 \setminus N(F)$ are incompressible.
\end{claim}

\begin{proof}
First, we claim that every component of $H(t_i)\cap F$ which is a loop is either essential
or inessential on both surfaces. If not, then as in the proof of
\fullref{Clm:UnlabelledInterval} there is a loop component $\delta$ of
$H(t_i)\cap F$ which bounds a compressing disk for $H(t_i)$. Assume that the compressing disk bounded by $\delta$ lies in $W(t_i)$, as the other case is similar. Pushing $\delta$ off of $F$ along $H(t_i)$ then yields a loop on $H(t_i) \setminus N(F)$ bounding a compressing disk in $W(t_i)$. This implies that there is a loop on $H(t_i-\epsilon) \setminus N(F)$ that bounds a compressing disk for $H(t_i-\epsilon)$ in $W(t_i-\epsilon)$. This violates the fact that the subinterval $(t_{i-1},t_i)$ does not have the label $\WW$.  

Now let $\Gamma_u$ denote the union of the inessential loops of $H(t_i) \cap F$ and $\Gamma _e$ the union of the essential loops. The intersection set $H(t_i) \cap F$ thus consists of $\Gamma _u$, $\Gamma _e$, and a figure eight curve $C$. Let $N_H(C)$ denote a closed neighborhood of $C$ on $H(t_i)$. If some component $\alpha$ of $\partial N_H(C)$ bounds a disk in $H(t_i)$ that contains $C$ then we say $C$ was {\it inessential}. 

Let $\pi\co H \times I \to H$ denote projection onto the first factor. Let $\pi_H=\pi \circ \Phi^{-1}$. Then, for each $t \in (0,1)$, the function  $\pi _H | H(t)$ is a map from $H(t)$ to $H$. 

The sets $\pi_H (\Gamma_u)$ and $\pi_H(\Gamma _e)$ are isotopic to subsets of $\pi _H (H(t_i-\epsilon) \cap F)$ and $\pi _H (H(t_i+\epsilon)\cap F)$, for sufficiently small $\epsilon$. Such an isotopy induces an identification of $\Gamma _u$ and $\Gamma _e$ with subsets of $H(t_i-\epsilon) \cap F$ and $H(t_i+\epsilon) \cap F$. Furthermore the loop $\alpha$ (if it exists) can be identified with loops on $H(t_i-\epsilon)$ and $H(t_i+\epsilon)$ which are disjoint from $F$. 

Let $H_0$, $H_-$ and $H_+$ denote the surfaces obtained by isotoping $H(t_i)$, $H(t_i-\epsilon)$ and $H(t_i+\epsilon)$, preserving $\Gamma _e$, but removing $\Gamma _u$. In each case these isotopies can be achieved via a series of identical moves on innermost disks. Note that if the figure eight $C$ is inessential and surrounded by some loop of $\Gamma_u$ then it will disappear in the course of these isotopies.

Now suppose $C$ was inessential but did not disappear (and is therefore not surrounded by some loop of $\Gamma _u$). By definition $\alpha$ bounds a disk $D$ on $H_0$ (which can be identified with disks on $H_-$ and $H_+$).  As $F$ is incompressible any intersection of $D$ with $F$ can be removed by a further isotopy of $H_0$, $H_-$ and $H_+$. Henceforth, we will assume that if $C$ is inessential then $\alpha$ bounds disks in $H_0$, $H_-$ and $H_+$ which are disjoint from $F$. 

Let $V_0, W_0, V_-, W_-, V_+$, and $W_+$ be the corresponding compression bodies bounded by $H_0$, $H_-$, and $H_+$. By assumption the interval $(t_{i-1},t_i)$ does not have the label $\WW$. It thus follows that no essential loop of $H_-$, disjoint from $F$, bounds a compressing disk in $W_-$. This is because only inessential loops are effected in the passage from $H(t_i-\epsilon)$ to $H_-$. Similarly we may conclude that no essential loop of $H_+$, disjoint from $F$, bounds a compressing disk in $V_+$. 

Assume, to obtain a contradiction, that $E'$ is a compressing disk for a component $H'$ of $H_0 \setminus N(F)$. Since every loop of $H_0 \cap F$ is essential on $H_0$, and $C$ was removed if it was inessential, it follows that $\partial E'$ is essential on $H_0$. Furthermore, as only the inessential intersection curves were effected in the passage from $H(t_i)$ to $H_0$ it follows that $\partial E'$ is an essential loop on $H(t_i)$, and is disjoint from $F$. It follows from \fullref{Lem:NoNesting} that there is a compressing disk $E$ for $H(t_i)$ with $\partial E=\partial E'$. Hence $\partial E$ is also disjoint from $F$. 

The loop $\partial E$ can be identified with essential loops of both $H(t_i-\epsilon) \setminus N(F)$ and $H(t_i+\epsilon) \setminus N(F)$ which bound similar compressing disks. We conclude the disk $E$ may be identified with a compressing disk for both $H_-$ and $H_+$ with boundary disjoint from $F$. If $E \subset W(t_i)$ then this violates the fact that there is no compressing disk for $H_-$ in $W_-$ with boundary disjoint from $F$. On the other hand, if $E \subset V(t_i)$, then we contradict the fact
that there is no compressing disk for $H_+$ in $V_+$ with boundary disjoint from $F$. 

We conclude that the components of $H_0 \setminus N(F)$ are
incompressible in the respective submanifolds of $M \setminus N(F)$, as
asserted by the third conclusion of the lemma.
\end{proof}

The third conclusion of \fullref{Lem:StronglyIrreducibleSweepout} follows. This completes the proof of \fullref{Lem:StronglyIrreducibleSweepout}.
\end{proof}

We now use the above result to establish the following lemma.

\begin{lemma}
\label{Lem:Sweepout}
Let $M$ be a compact, irreducible, orientable 3--manifold with $\partial M$ incompressible, if non-empty. Suppose $M=X \cup _F Y$, where $F$ is essential, connected, and closed. Suppose $M=V \cup _H W$, where $H$ is a Heegaard surface. Then
either $H$ is an amalgamation of splittings of $X$ and $Y$ or there
are properly embedded surfaces $H_X \subset X$ and $H_Y \subset Y$
with boundaries on $F$ such that at least one of the following holds:
\begin{enumerate}
\item The surfaces $H_X$ and $H_Y$ are incompressible, not boundary parallel, $\partial H_X=\partial H_Y$ and
  $\chi(H_X)+\chi(H_Y) \ge \chi(H)$.
\item After possibly exchanging $X$ and $Y$ the surface $H_X$ is
  incompressible, not boundary parallel, the surface $H_Y$ is
  strongly irreducible, $\partial H_X=\partial H_Y$ and
  $\chi(H_X)+\chi(H_Y) \ge \chi(H)$.
\item The surfaces $H_X$ and $H_Y$ are incompressible, not boundary parallel, $\partial H_X \cap \partial H_Y =\emptyset$,
  and $\chi(H_X)+\chi(H_Y)-1 \ge \chi(H)$.
\end{enumerate}
\end{lemma}

\begin{remark}
If $H$ is assumed to be strongly irreducible then we will show that
each of the above inequalities can be replaced by equalities.
\end{remark}

\begin{proof}
By Scharlemann--Thompson \cite{ScharlemannThompson94} we may {\it untelescope} the
Heegaard splitting $H$.  That is, there is a generalized Heegaard splitting 
$\{H_i\}_{i=0}^{2n}$ of $M$ with thick and thin levels obtained from $H$ by some number of compressions. Furthermore, we can find such a generalized Heegaard splitting such that each thick level $H_i$ is strongly irreducible in the submanifold of $M$ cobounded by $H_{i-1}$ and $H_{i+1}$. It is shown in \cite{ScharlemannThompson94} that in such a generalized Heegaard splitting each thin level is incompressible in $M$. 

Isotope $F$ to meet the set of thin levels of $\{H_i\}$ in a minimal
number of curves. Suppose first that for some $i$, the surface $F$ is
parallel to a component of the thin level $H_{2i}$. Then the
components of $\{H_i\}$ which meet $X$ form an untelescoped Heegaard
splitting of $X$, and the components which meet $Y$ form an
untelescoped Heegaard splitting of $Y$. Telescoping (the operation
which is the inverse of untelescoping) now produces Heegaard
splittings of $X$ and $Y$ with amalgamation $H$. Hence, the first
conclusion of \fullref{Lem:Sweepout} follows.

Now suppose $F$ intersects the thin level $H_{2i}$. Then $F$ divides
$H_{2i}$ into subsurfaces $H_X \subset X$ and $H_Y \subset Y$. We
claim that $H_X$ is incompressible in $X$ and $H_Y$ is incompressible
in $Y$. If not, then there is some compressing disk $D$ for $H_X$
(say) in $X$. As $H_{2i}$ is incompressible in $M$, $\partial D$
bounds a disk $E$ in $H_{2i}$. Since $\partial D$ is essential in $H_X$ but inessential in $H_{2i}$ the surface $F$ must intersect the disk $E \subset H_{2i}$. As $M$ is irreducible we can now do a
sequence of isotopies to remove all curves of $E \cap F$, reducing the
number of times $F$ meets the set of thin levels.

Since $F$ meets all thin levels minimally it also follows that neither
$H_X$ nor $H_Y$ are boundary parallel. Finally, since $H_{2i}=H_X \cup
H_Y$, and $H_{2i}$ is obtained from $H$ be some number of
compressions, we have $\chi(H_X)+\chi(H_Y) \ge \chi(H)$. Hence,
Case~(1) of the conclusion of \fullref{Lem:Sweepout} follows.

We are now reduced to the case where $F$ misses all thin levels, and
is parallel to none. Hence, $F$ is completely contained in a
submanifold with incompressible boundary which has a strongly
irreducible Heegaard splitting, obtained from $H$ by some number of
compressions. It suffices, then, to prove \fullref{Lem:Sweepout} in
the case where $H$ is strongly irreducible.

Use \fullref{Lem:StronglyIrreducibleSweepout} to isotope $H$ so that
it is transverse or almost transverse to $F$, and so that the conclusion of \fullref{Lem:StronglyIrreducibleSweepout} follows.  If $H$ is transverse to $F$ then let
$H_X=H \cap X$ and $H_Y=H \cap Y$, and Case~(1) or (2) of the lemma at hand follows. 

The remaining case
is when $H$ meets $F$ almost transversally. Let $p$ denote the saddle point of $H \cap F$. Isotope $H$ by pushing the point $p$
slightly into $Y$, to obtain the surface $H'$. Hence, $H'$ is transverse
to $F$. Furthermore, any compressing disk for $H_X=H' \cap X$ is a
compressing disk for $H \cap X$, so there must be none by \fullref{Lem:StronglyIrreducibleSweepout}. We conclude
$H_X$ is a properly embedded, incompressible surface in $X$.
Similarly, by pushing $p$ slightly into $X$ we may obtain from $H$ a
properly embedded, incompressible surface $H_Y \subset Y$.

        \begin{figure}[ht!]
        \labellist\small
        \pinlabel {$H$} at 357 223 
        \pinlabel {$F$} [b] at 262 408
        \pinlabel {$X$} at 45 437 
        \pinlabel {$Y$} at 36 335
        \pinlabel {$p$} [l] at 220 463
        \pinlabel {$H_X$} at 450 448 
        \pinlabel {$H_Y$} at 208 344
        \endlabellist
        \cl{\psfig{file=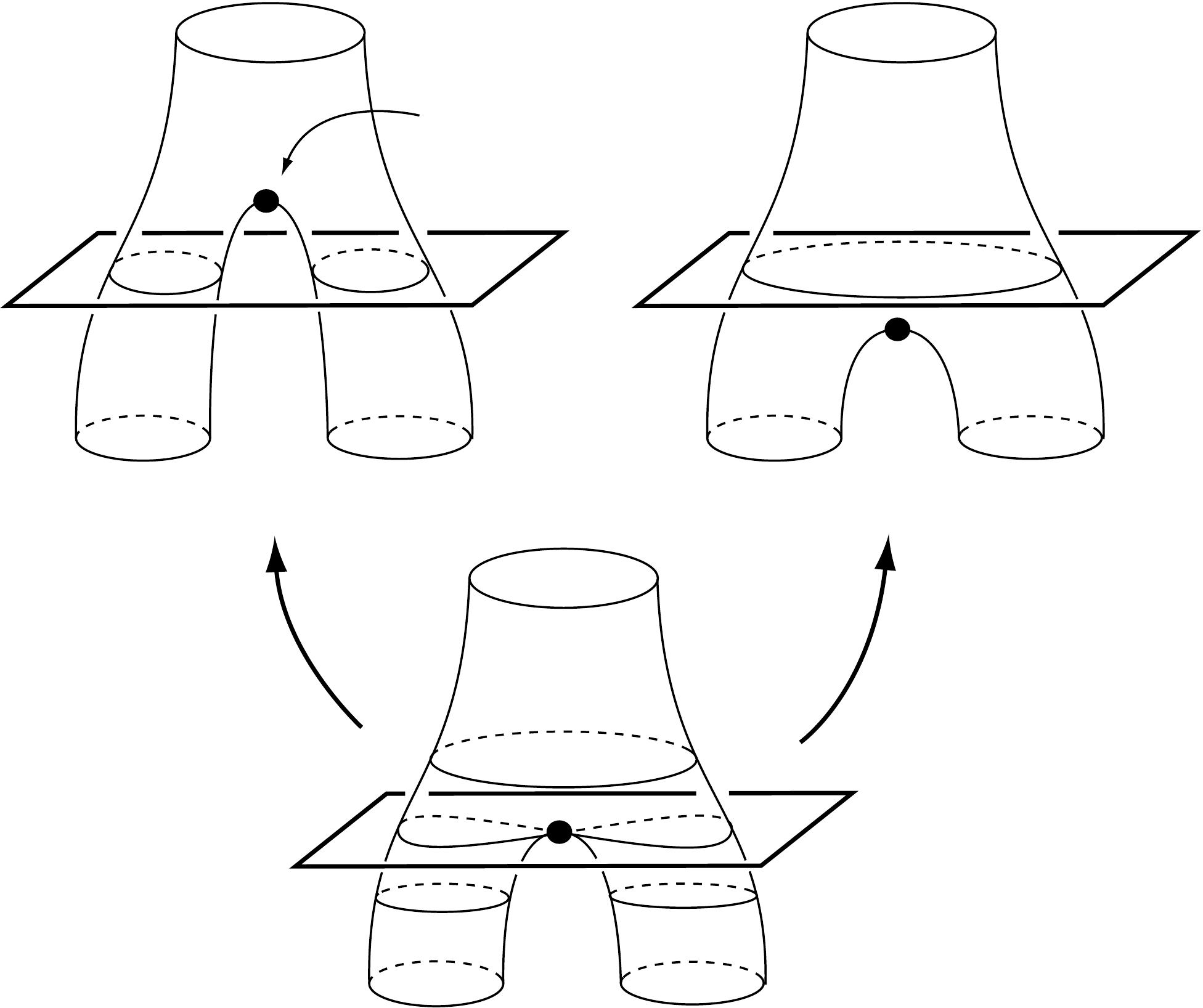,width=4in}}
        \caption{$H$ differs from $H_X \cup H_Y$ by a pair of pants.}
        \label{Fig:saddle}
        \end{figure}
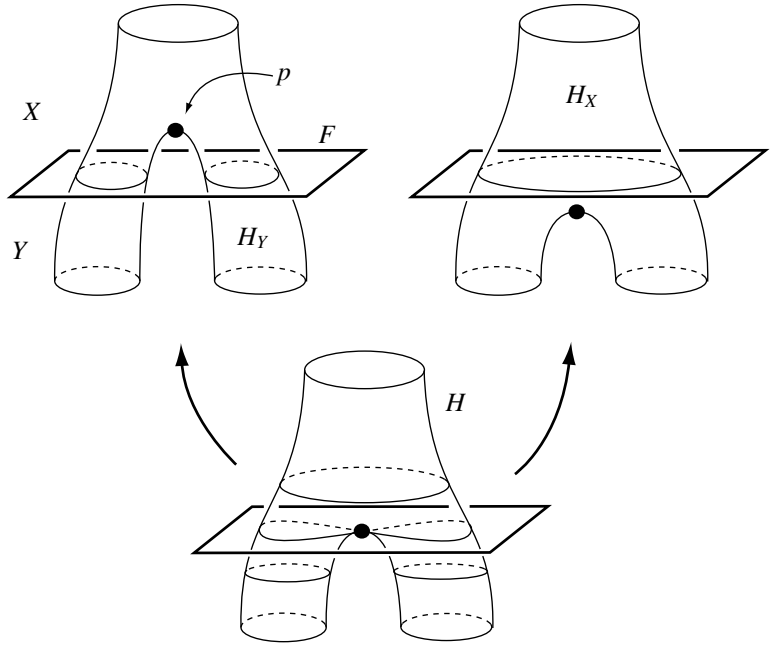

As $H$ and $F$ are orientable, it follows that $H_X \cap F$ may be
made disjoint from $H_Y \cap F$. Furthermore, the only essential
difference between $H_X \cup H_Y$ and $H$ is a pair of pants, having
Euler characteristic negative one (see \fullref{Fig:saddle}).  Hence, Case~(3) of the conclusion
of \fullref{Lem:Sweepout} now follows.
\end{proof}

\section{Manifolds with no strongly irreducible Heegaard splittings}
\label{s:saul}

A {\em knot manifold} is a compact, orientable, irreducible
three--manifold with a single boundary
component, which is incompressible and homeomorphic to a torus.  
The goal of this section is to prove the following theorem:

\begin{theorem}
\label{Thm:ManifoldsWithoutStronglyIrreducibleSplittings}
Suppose that $X$ and $Y$ are {\em knot manifolds} and $\varphi \from \bdy X
\to \bdy Y$ is a {\em sufficiently complicated} homeomorphism.  Then the manifold
$M(\varphi) = X \cup_\varphi Y$ has no strongly irreducible
Heegaard splittings.
\end{theorem}

Note the similarity of \fullref{Thm:ManifoldsWithoutStronglyIrreducibleSplittings} to Cooper and Scharlemann's result \cite{cs:99}. That paper proves that if a 3--manifold is constructed by identifying the boundary components of $T^2 \times I$ via a ``sufficiently complicated" map then there are no strongly irreducible  Heegaard splitting of the resulting 3--manifold.

To make the statement of \fullref{Thm:ManifoldsWithoutStronglyIrreducibleSplittings}  precise we must give a reasonable definition of the term  {\em sufficiently complicated}. To this end fix, once and for all,  psuedo-triangulations of $X$ and $Y$ with one vertex. (A {\it psuedo-triangulation} is a decomposition into simplices where any two such simplices intersect in a collection of lower dimensional simplices.)  Let $\Delta(X)$ be the set of slopes in $\partial X$ which are the boundary of some normal or almost normal surface in $X$.  Note that $\Delta(X)$ is finite, by a result of Jaco and
Sedgwick~\cite{JacoSedgwick03} (see also Theorem 9.7 of Bachman
\cite{bachman:01} for a discussion of the almost normal case). Define $\Delta(Y)$ similarly.

Recall now the definition of the {\em Farey graph}, $\calF(X)$.  The
vertices of $\calF(X)$ are all slopes in $\bdy X$.  Two slopes are connected by an edge if they intersect once. The {\it distance} between two slopes is then defined to be the minimal number of edges required in a path connecting them. The {\it distance} between two sets of slopes is the minimal distance between their elements. 

\begin{define}
A map $\varphi \from \bdy X \to \bdy Y$ is {\em sufficiently complicated} if the distance between $\Delta(X)$ to $\varphi^{-1}(\Delta(Y))$ inside of $\calF(X)$ is at least two.
\end{define}

\begin{remark}
\label{Rem:Ubiquity}
Note that, as $\Delta(X)$ and $\Delta(Y)$ are finite, ``most''
elements of $\MCG(\TT^2)$ $\isom SL(2, \ZZ)$ are sufficiently complicated,
in the above sense.  In particular any sufficiently large power of an
Anosov map is sufficiently complicated.  The same holds for all but a
finite number of Dehn twists.  
\end{remark}

Before giving the proof of
\fullref{Thm:ManifoldsWithoutStronglyIrreducibleSplittings} we
must discuss boundary compressions.  Suppose $G \subset N$ is a
properly embedded, two--sided surface in a compact, orientable,
irreducible three--manifold $N$.  We suppose further that $\bdy N$ is
incompressible in $N$.  Suppose $D \subset N$ is a boundary
compression for $G$.  

\begin{define}
The boundary compression $D$ is {\em honest} if $D \cap \partial N$ is
essential as a properly embedded arc in $\bdy N \setminus \bdy G$.  If
$D$ is not honest it is {\em dishonest}.
\end{define}

\begin{define}
Let $N$ be a knot manifold.  We now define the {\em
banding}, $\band{D}$, of a boundary compression $D$ for $G$. First assume $D$ is honest. Then $D \cap \partial N$ meets distinct boundary components of $\bdy G$,
as $G$ is orientable.  These components of $\bdy G$ cobound an annulus
$A \subset \bdy N$ such that $D \cap \partial N\subset A$.  Let $D'$ denote the disk obtained from $A$ by removing a neighborhood of $D \cap \partial N$ and attaching two parallel copies of $D$.  Isotope $D'$ to be disjoint
from $\bdy N$ while maintaining $\bdy D' \subset G$.  The resulting
disk is the desired banding $\band{D}$ of $D$.

Now suppose $D$ is dishonest. Then the arc $D \cap \partial N$ cobounds, with a subarc of $\partial G$, a subdisk $D'$ of $\partial N$. The disk $\band{D}$ is obtained by pushing $D''=D \cup D'$ into the interior of $N$, while maintaining $\partial D'' \subset G$. 
\end{define}

Note that when $C$ is a compressing disk and $D$ is a boundary-compressing disk (honest or dishonest) if $C \cap D=\emptyset$ then $C \cap \band{D}=\emptyset$.

\begin{lemma}
\label{Lem:Band}
If $D$ is a boundary compression for $G$ and $\bdy N = \TT^2$ then $G$
is either compressible or the component of $G$ meeting $D$ is a boundary parallel annulus.\hfill$\sq$
\end{lemma}

Recall that by a {\em strongly irreducible} surface we mean a
properly embedded, two--sided surface which compresses on both sides and all pairs of
compressing disks on opposite sides must meet. We now strengthen this definition to account for boundary compressions, as in Bachman \cite{bachman:01}.

\begin{define}
A properly embedded, separating surface is {\it $\partial$--strongly irreducible} if 
\begin{enumerate}
	\item every compressing and boundary-compressing disk on one side meets every compressing and boundary-compressing disk on the other side, and
	\item there is at least one compressing or boundary-compressing disk on each side.
\end{enumerate}
\end{define}

\begin{lemma}
\label{Lem:StronglyIrreducibleImpliesDistanceOne}
Let $N$ be a knot manifold. Let $G$ be a separating, properly embedded, connected surface in $N$ which is strongly irreducible, has non-empty boundary, and is not peripheral. Then either $G$ is $\partial$--strongly irreducible or $\partial G$ is at most distance one from the boundary of some properly embedded surface which is both incompressible and boundary-incompressible.
\end{lemma}

\begin{proof}
Suppose $G$ divides $N$ into $V$ and $W$. If $G$ is not $\partial$--strongly irreducible then there are disjoint disks $D \subset V$ and $E \subset W$ such that at least one, say $D$, is a boundary-compressing disk. The disk $E$ is either a compression or a boundary compression. 

Since $G$ is not a boundary parallel annulus we know by \fullref{Lem:Band} that the banding disk $\band{D}$ is a compressing disk for $G$. If $E$ is a compressing disk then $E \cap D=\emptyset$ implies that $E \cap \band{D}=\emptyset$, contradicting strong irreducibility. We conclude $E$ is a boundary compression. 

Let $G'$ denote the result of boundary-compressing $G$ along both $D$ and $E$. Let $V'$ and $W'$ denote the sides of $G'$ which correspond to $V$ and $W$. We now claim that $G'$ is incompressible. Suppose $D'$ is a compressing disk for $G'$ in $V'$. Then $D'$ must have been a compressing disk for $G$ in $V$ which was disjoint from $E$, and hence disjoint from $\band{E}$. This contradicts the strong irreducibility of $G$. By symmetry we conclude $G'$ is incompressible. 

We now claim $G'$ is boundary incompressible as well. Suppose $C$ is a boundary-compressing disk for $G'$. Since $G'$ is incompressible we know $\band{C}$ is not a compressing disk, so it follows from \fullref{Lem:Band} that $G'$ must be a boundary parallel annulus. It follows that all of $G$ was isotopic into a neighborhood of $\partial N$, contradicting our hypotheses. 

It remains only to show that $\partial G$ is at a distance of at most one from $\partial G'$. In order for the slope of $\partial G'$ to be different from the slope of $\partial G$ all of the loops of $\partial G$ must meet either $D$ or $E$. This immediately implies $|\partial G| \le 4$. The possibility that $|\partial G|$ is one or three is ruled out by the fact that $G$ is separating. The fact that $D$ and $E$ are on opposite sides of $G$ rules out $|\partial G|=4$, since we are assuming that every component of $\partial G$ meets either $D$ or $E$.  

If $|\partial G|=2$, both $D$ and $E$ are dishonest, and each meets different components of $\partial G$ then $\band{D} \cap \band{E}=\emptyset$. This violates the strong irreducibility of $G$. 

        \begin{figure}[ht!]
        \vspace{0 in}
        \cl{\includegraphics[width=3in]{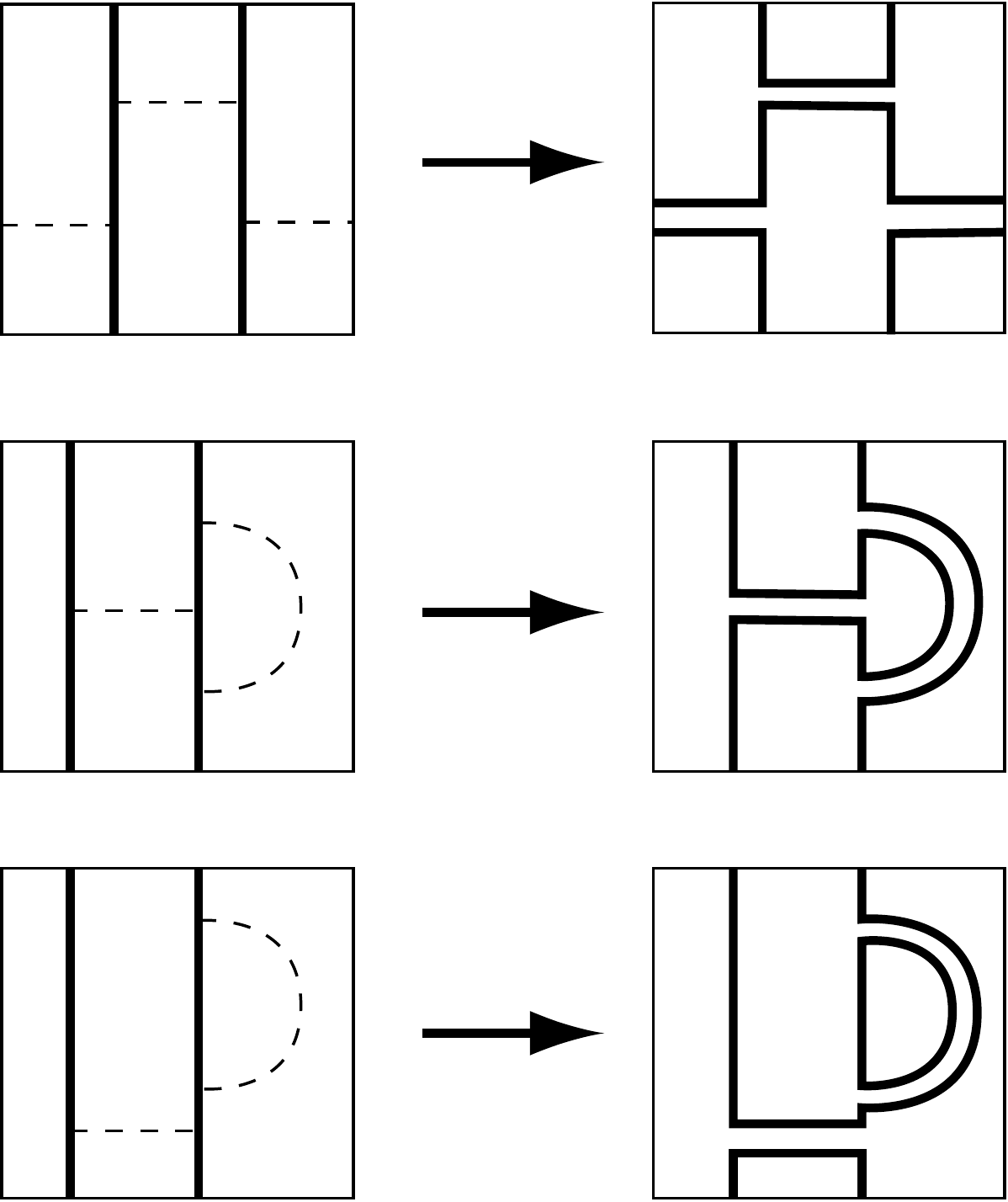}}
        \caption{Possible effects of boundary-compression on $\partial G$}
        \label{Fig:twobdy}
        \end{figure}

There are three remaining cases. In each of these cases $|\partial G|=2$ and both boundary loops are affected by the transition to $G'$. See \fullref{Fig:twobdy}. In the top picture both $D$ and $E$ are honest. The two loops of $\partial G$ are transformed into two loops, both distance one from the original. In the middle picture exactly one of the disks $D$ or $E$ is dishonest, and the boundary slope remains unchanged. The configuration depicted at the bottom of \fullref{Fig:twobdy} cannot happen, since it represents a situation in which $\band{D}$ is disjoint from $\band{E}$, contradicting the strong irreducibility of $G$. 
\end{proof}

We conclude with:

\begin{proof}[Proof of
\fullref{Thm:ManifoldsWithoutStronglyIrreducibleSplittings}]
Suppose that $X$ and $Y$ are triangulated knot manifolds, as above.
Fix a gluing $\varphi \from \bdy X \to \bdy Y$.  Suppose that $H
\subset M(\varphi) = X \cup_\varphi Y$ is a strongly irreducible
Heegaard splitting surface.  Let $F \homeo \TT^2$ be the image of
$\bdy X$ inside of $M(\varphi)$.

Now apply \fullref{Lem:StronglyIrreducibleSweepout} and
Remark~\ref{Rmk:HTransverseToF} to the pair $H$ and $F$ in
$M(\varphi)$. Let $H_X$ be a component of $H \cap X$ which is incompressible and not
a boundary parallel annulus, if such exists.  If no such component
exists take $H_X$ to be the non-boundary parallel component of $H \cap
X$.  In this case $H_X$ is strongly irreducible.  (At least one component
of $H \cap X$ is not boundary parallel. Otherwise $H$ is isotopic into $Y$, 
a contradiction.)  Choose $H_Y$ similarly and note that,
by \fullref{Lem:StronglyIrreducibleSweepout}, not both of $H_X$ and
$H_Y$ are strongly irreducible.  Note that $\bdy H_X$ and $\varphi^{-1}(\bdy H_Y)$ have the same slope. 

Suppose that $H_X$ and $H_Y$ are both incompressible.  As $\bdy X
\homeo \bdy Y \homeo \TT^2$ it follows from \fullref{Lem:Band} that
$H_X$ and $H_Y$ are also boundary incompressible.  So $H_X$ and $H_Y$
may be normalized with respect to the given
triangulations Haken~\cite{haken:61}.  It follows that the sets $\Delta(X)$ and
$\varphi^{-1}(\Delta(Y))$ intersect and thus $\varphi$ is not sufficiently
complicated.

Suppose now that $H_X$ is incompressible and thus
boundary incompressible.  Suppose that $H_Y$ is a strongly irreducible
surface.  Then, by \fullref{Lem:StronglyIrreducibleImpliesDistanceOne}, either $H_Y$ is $\partial$--strongly irreducible or $\bdy H_Y$ intersects the boundary of some incompressible, boundary incompressible surface $H'_Y$ at most once. In the latter case $H_Y'$ may be normalized, and hence $\partial H'_Y \in \Delta(Y)$. In the former case it follows from work of the first author (Corollary 8.9 of \cite{bachman:01}) that the surface $H_Y$ is
properly isotopic to an almost normal surface, and so $\partial H_Y \in \Delta(Y)$. In either case we see $\partial H_X$ (an element of $\Delta(X)$) is within distance one from some element of $\varphi^{-1}(\Delta (Y))$ and hence $\varphi$ is not sufficiently complicated.
\end{proof}

\section{Amalgamating small manifolds}
\label{s:eric}

Let $X$ be a manifold with boundary.  The {\it tunnel number } of
$X$, $t(X)$ is the minimal number of properly embedded arcs that
need to be drilled out of $X$ to obtain a handlebody; i.e.\ so
that $X  \setminus  N({\rm arcs})$ is a handlebody.  The {\it handle
number} of $X$ is the minimal number of properly embedded arcs
that need to be drilled out of $X$ to obtain a compression body;
i.e.\ so that $X \setminus N({\rm arcs})$ is a compression body. If $|\partial X|=1$ then $t(X)=h(X)$.

Let $M = X \cup_F Y$ be a manifold obtained by gluing $X$
and $Y$, two connected small manifolds with incompressible boundary,
along a collection of boundary components homeomorphic to a surface
$F$.   The goal of this section is to show that the Heegaard genera
of $X$ and $Y$ are bounded in terms of the Heegaard genus of $M = X
\cup_F Y$. More specifically, we establish:

\begin{theorem}
\label{Cor:genusbounds}
Let $M$ be a compact, orientable 3--manifold with incompressible\break
boundary.  Suppose $M$ is obtained by gluing two connected, small manifolds
along a union of incompressible boundary components, $M = X \cup_F Y$.  Then the
following statements hold:
\begin{enumerate}
	\item $g(M) \geq \frac{1}{2}(h(X) + h(Y)) $ 
	\item if $M$ is closed and $F$ is connected, $g(M) \geq \frac{1}{2}(t(X) + t(Y))$ 
	\item $g(M) \geq \frac{1}{2} ( g(X) + g(Y) -
2g(F))$.
\end{enumerate}
\end{theorem}

The theorem is motivated by the fact that a properly embedded,  incompressible surface cuts a small manifold into one or two
compression bodies. 

We begin with the following definitions. Let $F$ be an orientable surface, possibly
with boundary components, and possibly disconnected.  Let $C$ be
the manifold obtained by forming $F \times I$ and attaching one
handles to the surface $F \times \{1\}$. Then $C$ is a {\it
relative compression body}.  We label the boundary as follows: the
{\it negative boundary} is $\partial_-C = F \times \{0\}$, the {\it vertical
boundary} is $\partial_V C= \partial F \times I$, and the {\it positive
boundary} is $\partial_+ C =\overline{\partial C  \setminus (\partial_-C \cup
\partial_V C)}$.  The vertical boundary is a collection of annuli.
It is important to note that a given manifold may admit many
relative compression body structures. For example, if $F$ is a
surface with boundary and $C=F \times I$, then $C$ can be thought
of as a relative compression body with $\partial_- C = F \times
\{0\}$, or $C$ can be thought of as a handlebody with $\partial_-
C = \emptyset$.  In fact, given a relative compression body $C$,
it is always possible to think of $C$ as a (non-relative)
compression body by {\it promoting} all non-closed components of
$\partial_- C$ and all components of $\partial_V C$ to the positive
boundary.

A {\it relative Heegaard splitting} is the union of two
relative compression bodies, identified along their positive
boundaries.  The splitting will be considered {\it non-trivial} if
neither relative compression body is a product; i.e.\ both
compression bodies have 1--handles.

\begin{lemma}
\label{Lem:SI} Let $X$ be a manifold that admits a non-trivial, strongly irreducible and 
relative Heegaard splitting $X = C_1 \cup C_2$. Then $\partial_-C_1$ and
$\partial_-C_2$ are incompressible in $X$.
\end{lemma}

\begin{proof}
An examination of the proof of the Haken Lemma \cite{haken:68} (see also Jaco
\cite{jaco:80})
will reveal that it applies directly to the case of relative
Heegaard splittings.  In particular, if either $\partial_-C_1$ or
$\partial_-C_2$ has compressible boundary, then there is a
compressing disk $D$ for the boundary component that meets the
splitting surface in a single closed loop.  The loop decomposes
the compressing disk into a vertical annulus in one compression
body, say $C_1$, and a disk $D_2 \subset C_2$. Since $C_1$ is not
a product we can find a compressing disk $D_1$ for $\partial_+
C_1$, disjoint from the annulus, and hence disjoint from $D_2$.
The pair $(D_1,D_2)$ contradicts strong irreducibility of the relative Heegaard splitting.
\end{proof}

\begin{lemma}
\label{Lem:CompBody} An irreducible connected small manifold with
compressible boundary is a compression body.
\end{lemma}

\begin{proof}
Let $X$ be a connected small manifold with compressible boundary.
In an optimistic fashion, denote a compressible boundary component
by $\partial_+ X$ and all other components by $\partial_- X$.
Since $\partial_+ X$ is compressible it bounds a (not properly embedded) submanifold $C$ of $X$ which is a compression body, so that $\partial_+C =
\partial_+X$. Choose $C$ to be maximal in this regard. Precisely,
choose $C$ so that $\partial_-C$ contains no 2--spheres ($X$ is
irreducible) and so that $\sum (1-\chi(S_i))$ is minimal, where $\{S_i\}$ are the components of $\partial _-C$.

If $S$ is a component of $\partial_-C$ then  $S$ is incompressible in
$C$.  Suppose $D$ is a compressing disk for $S$ in $X \setminus C$. Then $D$ is the core of a 2--handle that we can attach to $C$ to obtain a new compression body with negative boundary ``smaller" than that of $C$. This contradicts our minimality assumption. We conclude $S$ is incompressible in $X \setminus C$. As $X$ is small, $S$ must be peripheral, and
since $C$ is not a product, it is parallel in $X \setminus C$ to a
component of $\partial_-X$. The (possibly disconnected) surface
$\partial_-C$ separates the components of $\partial_-X$ from
$\partial_+X$, so each component of $\partial_-X$ is in
fact parallel to a component of $\partial_-C$. The parallelism
yields an isotopy between $X$ and $C$.  $X$ is therefore a
compression body.  Note that only one boundary component,
$\partial_+X$, is compressible.
\end{proof}

\begin{theorem}
\label{Prop:Incomp}
Let $H_X$ be a non-peripheral, connected, incompressible surface that is properly
embedded in a connected, small manifold $X$. Then $h(X) \leq 1 - \chi(H_X)$. If $X$
has a single boundary component or $H_X$ meets every boundary component
of $X$, then this applies to the tunnel number: $t(X) \leq 1 -
\chi(H_X)$.
\end{theorem}

\begin{proof}
Let $\partial_1 X$ denote those boundary components of $X$ that meet
$H_X$ and $\partial_2 X$ denote those boundary components which do not
meet $H_X$.

        \begin{figure}[ht!]
        \labellist\small
        \pinlabel {$\partial _1 X$} [r] at 156 426
        \pinlabel {$\partial _2 X$} [tl] at 334 429
        \pinlabel $X_1$ at 431 513 
        \pinlabel $X_2$ at 428 460 
        \pinlabel $H_X$ [bl] at 311 622
        \pinlabel $H_X'$ [tr] at 178 370
        \endlabellist
        \cl{\hspace{.7in}\psfig{file=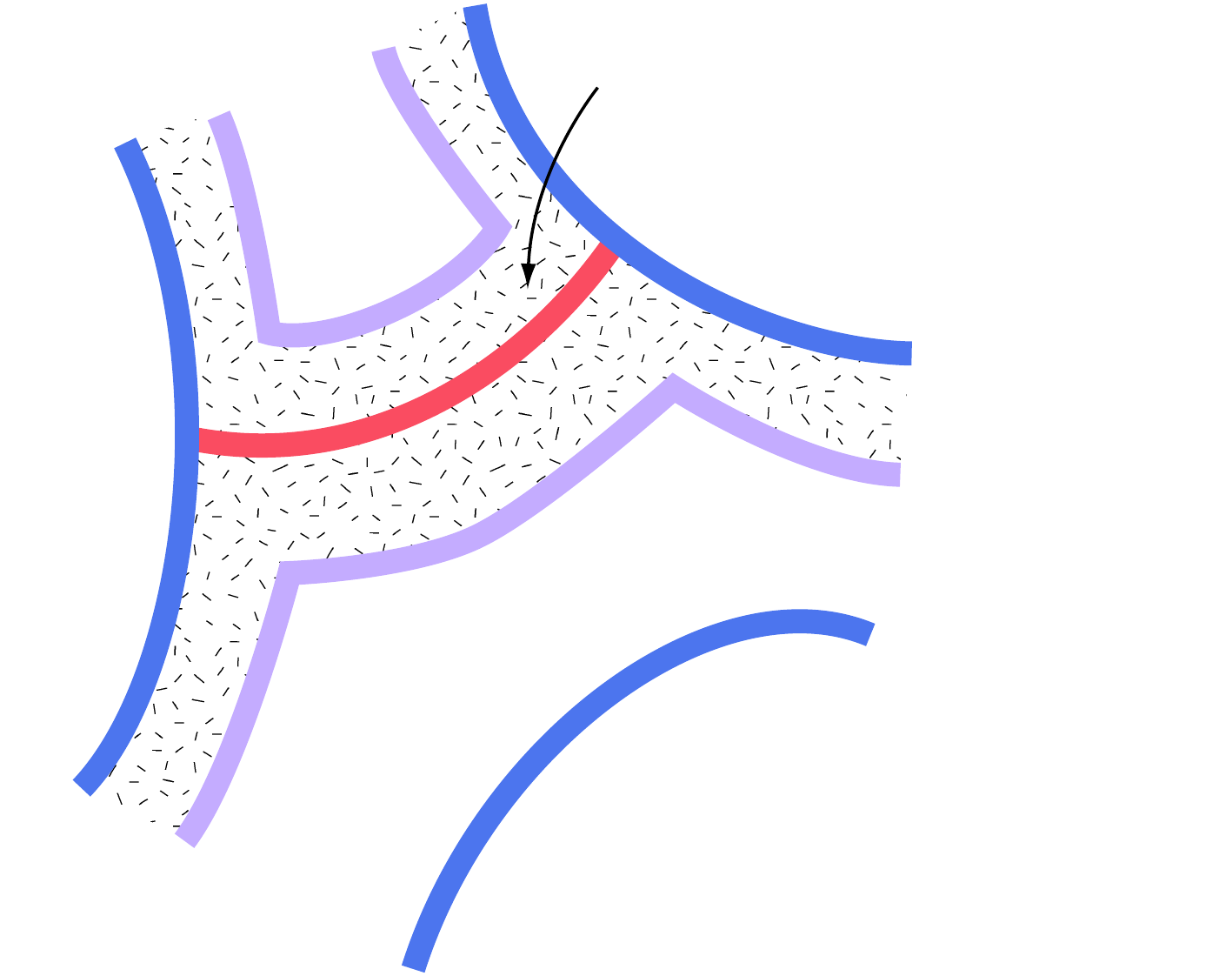,width=3in}}
        \caption{Labelling the boundary components of $X$}
        \label{Fig:decomp}
        \end{figure}
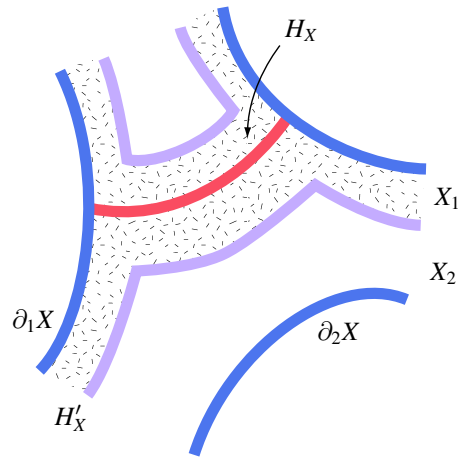

Let $X_1 = \overline{N(H_X \cup \partial_1 X)}$ and $X_2 = \overline {X
\setminus X_1}$. This decomposes $X$ into $X=X_1 \cup_{H_X'}
X_2$, where $H_X'$ is the common boundary
of $X_1$ and $X_2$.  See the schematic in Figure~\ref{Fig:decomp}.
Note that $\partial_1X$ and $H_X$ are contained in $X_1$ and
$\partial_2X$ is contained in $X_2$. Since $H_X$ is connected it follows that $X_1$ is 
connected. If $H_X$ separates $X$ then $X_2$ will have two components.

The surface $H_X'$ will have two components if $H_X$ separates and one
component otherwise.  Since $X$ is a small manifold, each component of
$H_X'$ is either compressible in $X$ or peripheral to a boundary
component of $X$.

\begin{proofclaim}
{\it If a component of $H_X'$ is compressible, it is compressible into $X_2$.}
\end{proofclaim}

\begin{proof}
If there is a compressing disk for the compressible component of $H_X'$ then there is one that is disjoint from $H_X$. This is because any intersection could be removed by surgery. If $H_X'$ has  two components then $H_X$ separates them. Hence our chosen compressing disk does not meet the other component of $H_X'$.   Therefore, our compressing disk is
properly embedded in either $X_1\setminus N(H_X)$ or $X_2$.  But, $X_1\setminus N(H_X)$ is
a product and has incompressible boundary.  It follows that a
compressible component of $H_X'$ is compressible into $X_2$.
\end{proof}

\begin{proofclaim}
{\it No component of $H_X'$ is peripheral into $\partial_1 X$.}
\end{proofclaim}

\begin{proof}
If this occurred, $X_1$ would be contained in a product neighborhood
of a boundary component.  This in turn implies that $H_X$ was
peripheral.
\end{proof}

\begin{proofclaim}
{\it Each component of $X_2$ is a compression body.}
\end{proofclaim}

\begin{proof}
Suppose that a component $X'$ of $X_2$ contains a closed non-peripheral essential
surface $G$. Since $X$ is small, $G$ is either compressible in $X$
or parallel to a component $G'$ of $\partial X \setminus \partial X'$. In the latter case $G' \subset \partial_1 X$ or $G' \subset \partial _2 X \setminus \partial X'$. If $G' \subset \partial _2 X \setminus \partial X'$ then $H_X$ separates $G$ from $G'$. 

 Since $H_X$
is incompressible, any compressing disk $D \subset X$ for $G$ can
be isotoped so that it does not intersect $H_X$, and so can be isotoped to miss $X_1$.
Therefore $G$ is compressible in $X_2$, contradicting the essentiality of $G$. If $G' \subset \partial_1 X$ or $G' \subset \partial _2 X \setminus \partial X'$  then there
is a product containing $H_X$.  In particular, this implies that $H_X$
is contained in a product neighborhood of $\partial X$, contradicting the
fact that $H_X$ is not peripheral.  Thus, $X_2$ is small.

Each component of $H_X'$ is therefore compressible into $X_2$ or
parallel to a component of $\partial_2X$.  In either case, by
\fullref{Lem:CompBody} or by parallelism, $H_X' = \partial_+X_2$,
where $X_2$ is either one or two compression bodies.
\end{proof}

It is now straightforward to build a handle system for $X$ (tunnel
system in the case that $\partial_1X = \partial X$).  Choose
$\tau$, a minimal collection of arcs that are properly embedded in
$H_X$ and that cut $H_X$ into a single disk $D$.  The collection
$\tau$ contains $1-\chi(H_X)$ arcs.  Moreover, $\tau$ is a handle
system that induces a Heegaard splitting, $X=C_1 \cup C_2$, where
$C_1 = \overline{N(\partial_1 X \cup \tau)}$ and $C_2 =
\overline{X \setminus C_1}$. Clearly $C_1$ is a compression body.
$C_2$ is a compression body because it is formed by attaching a
1--handle (a neighborhood of the cocore of $D$) to the positive
boundary of the compression body/bodies $X_2$. This completes the proof of \fullref{Prop:Incomp}.
\end{proof}

\begin{theorem}
\label{Prop:SI}
Let $H_X$ be a non-peripheral, bi-compressible, connected, strongly irreducible
surface properly embedded in a connected, small manifold $X$. Then $h(X)
\leq 1 - \chi(H_X)$.  If $X$ has a single boundary component, then this
applies to the tunnel number: $t(X) \leq 1 - \chi(H_X)$.
\end{theorem}

\begin{proof}
We may apply the previous theorem if $X$ also contains a
non-peripheral incompressible surface with boundary whose negative Euler characteristic is less than that of $H_X$. We may therefore assume that $H_X$ is
a separating surface; if not we may compress $H_X$ to obtain such an
incompressible surface. As before we will let $\partial_1 X$
denote those boundary components of $X$ that meet $H_X$ and
$\partial_2 X$ denote those boundary components which do not meet
$H_X$.

By compressing $H_X$ maximally to both sides, we define a relative
Heegaard splitting of a submanifold $X' = C_1 \cup_{H_X} C_2 \subset X$.
Since we have compressed maximally, the negative boundary
components of $C_1$ and $C_2$ are incompressible outside $X'$.
They are incompressible inside $X'$ by \fullref{Lem:SI}. If any
component is non-peripheral, we have our conclusion via
\fullref{Prop:Incomp}.  Each component of $\partial_- C_i,
i=1,2$ is therefore peripheral. It now follows from the fact that $H_X$ is non-peripheral that
$X'$ is isotopic to $X$.

As in the earlier theorem, this structure defines a handle system
for $X$.  Choose $\tau$, a minimal collection of arcs that are
properly embedded in $H_X$ and that cut $H_X$ into a single disk $D$.
Now, $\tau$ is a handle system for $X$ that induces the Heegaard
splitting, $X=C_1' \cup C_2'$, where $C_1' =
\overline{N(\partial_1 X \cup \tau)}$ and $C_2' = \overline{X
\setminus C_1'}$. Clearly $C_1'$ is a compression body. $C_2'$ is a
compression body because it can be obtained by first promoting the
vertical and non-closed negative boundary components of $C_1$ and
$C_2$ and then joining the positive boundary of these
(non-relative) compression bodies with a $1$--handle (a
neighborhood of the cocore of $D$).

The handle number of $X$ is thus bounded by $1 -\chi(H_X)$.
\end{proof}

\begin{proof}[Proof of \fullref{Cor:genusbounds}]
Let $H$ be a minimal genus splitting of $M$.  If $H$ is an
amalgamation of splittings of $X$ and $Y$, then the result holds
trivially.  Otherwise, by \fullref{Lem:Sweepout} we can
construct properly embedded non-boundary parallel surfaces $H'_X
\subset X$ and $H'_Y \subset Y$ so that each is either
incompressible or strongly irreducible. As neither surface is boundary-parallel they contain components $H_X \subset H'_X$ and $H_Y \subset H'_Y$ which are non-boundary parallel and either incompressible or strongly irreducible. Furthermore, $\chi(H_X) +
\chi(H_Y) \geq \chi(H) = 2 - 2g(M)$, or equivalently, $g(M) \geq
\frac{1}{2}(2-\chi(H_X)-\chi(H_Y))$.

By either \fullref{Prop:Incomp} or \fullref{Prop:SI},
$X$ and $Y$ admit handle systems that are attached to components of
$F$ and so that the number of handles is at most $1 - \chi(H_X)$ and
$1 - \chi(H_Y)$, respectively.  The first two assertions of \fullref{Cor:genusbounds} follow.

Our induced splitting of $X$ is obtained by attaching $1-\chi(H_X)$ handles to $F$. The genus of $X$ is therefore bounded by 
\[g(X) \le g(F)+1-\chi(H_X).\]
Since a symmetric bound holds for $g(Y)$ we obtained the third conclusion of \fullref{Cor:genusbounds}.
\end{proof}

\bibliographystyle{gtart} \bibliography{link}

\end{document}